\documentclass{amsart}
\usepackage{amsthm}
\usepackage{amssymb}
\usepackage{amsmath}
\usepackage{graphicx,color}
\usepackage{hyperref}

\theoremstyle{plain}
\newtheorem{theorem}{Theorem}

\newtheorem{fact}{Fact}
\newtheorem{corollary}{Corollary}

\theoremstyle{definition}
\newtheorem{definition}{Definition}

\newtheorem{notation}{Notation}

\theoremstyle{remark}

\newcommand{\sub}{\operatorname{Sub}}
\newcommand{\sign}{\operatorname{sign}}

\newcommand{\Two}{
\,\begin{picture}(28,8)
\put(5,3){\circle{10}}
\put(25,3){\circle{10}}
\put(8.5,0){\vector(1,0){12}}
\put(20.5,6){\vector(-1,0){12}}
\end{picture}~
}

\newcommand{\foneA}{
\,\begin{picture}(28,8)
\put(5,3){\circle{10}}
\put(25,3){\circle{10}}
\put(8.5,0){\vector(1,0){12}}
\put(20.5,6){\vector(-1,0){12}}
\put(10,2.5){\vector(-1,0){10}}
\end{picture}~
}

\newcommand{\foneB}{
\,\begin{picture}(28,8)
\put(5,3){\circle{10}}
\put(25,3){\circle{10}}
\put(8.5,6){\vector(1,0){12}}
\put(20.5,0){\vector(-1,0){12}}
\put(10,2.5){\vector(-1,0){10}}
\end{picture}~
}
\newcommand{\foneC}{
\,\begin{picture}(28,8)
\put(5,3){\circle{10}}
\put(25,3){\circle{10}}
\put(9,0){\vector(3,1){11}}
\put(9.5,3){\vector(4,-1){11}}
\put(20.5,6){\vector(-1,0){12}}
\end{picture}~
}
\newcommand{\foneCII}{
\,\begin{picture}(28,8)
\put(5,3){\circle{10}}
\put(25,3){\circle{10}}
\put(9,0){\vector(3,1){11}}
\put(9.5,3){\vector(4,-1){11}}
\put(20.5,6){\vector(-1,0){12}}
\put(13,8){\tiny $\epsilon$}
\put(21,2){\tiny $+$}
\put(21,-1){\tiny $-$}
\end{picture}~
}

\newcommand{\foneD}{
\,\begin{picture}(28,8)
\put(5,3){\circle{10}}
\put(25,3){\circle{10}}
\put(8.5,6){\vector(1,0){12}}
\put(9,0){\vector(3,1){11}}
\put(9.5,3){\vector(4,-1){11}}
\end{picture}~
}

\newcommand{\foneDII}{
\,\begin{picture}(28,8)
\put(5,3){\circle{10}}
\put(25,3){\circle{10}}
\put(8.5,6){\vector(1,0){12}}
\put(9,0){\vector(3,1){11}}
\put(9.5,3){\vector(4,-1){11}}
\put(13,8){\tiny $\epsilon$}
\put(21,2){\tiny $+$}
\put(21,-1){\tiny $-$}
\end{picture}~
}

\newcommand{\One}{
\,\begin{picture}(28,8)
\put(5,3){\circle{10}}
\put(25,3){\circle{10}}
\put(20.5,3){\vector(-1,0){11}}
\end{picture}~
}

\newcommand{\Oneep}{
\,\begin{picture}(28,8)
\put(5,3){\circle{10}}
\put(25,3){\circle{10}}
\put(20.5,3){\vector(-1,0){11}}
\put(13,5){\tiny $\epsilon$}
\end{picture}~
}

\newcommand{\Oneb}{
\,\begin{picture}(28,8)
\put(5,3){\circle{10}}
\put(25,3){\circle{10}}
\put(10,3){\vector(1,0){10}}
\end{picture}~
}

\newcommand{\Onebep}{
\,\begin{picture}(28,8)
\put(5,3){\circle{10}}
\put(25,3){\circle{10}}
\put(10,3){\vector(1,0){10}}
\put(13,5){\tiny $\epsilon$}
\end{picture}~
}

\newcommand{\fone}{
\,\begin{picture}(28,8)
\put(5,3){\circle{10}}
\put(25,3){\circle{10}}
\put(8.5,0){\vector(1,0){12}}
\put(20.5,6){\vector(-1,0){12}}
\put(10,2.5){\vector(-1,0){10}}
\end{picture}~
}

\newcommand{\fonep}{
\,\begin{picture}(28,8)
\put(5,3){\circle{10}}
\put(25,3){\circle{10}}
\put(8.5,0){\vector(1,0){12}}
\put(20.5,6){\vector(-1,0){12}}
\put(10,2.5){\vector(-1,0){10}}
\put(2,3.8){\tiny $+$}
\end{picture}~
}
\newcommand{\fonepEP}{
\,\begin{picture}(28,8)
\put(5,3){\circle{10}}
\put(25,3){\circle{10}}
\put(8.5,0){\vector(1,0){12}}
\put(20.5,6){\vector(-1,0){12}}
\put(10,2.5){\vector(-1,0){10}}
\put(2.5,4){\tiny $\epsilon$}
\end{picture}~
}

\newcommand{\fonepPP}{
\,\begin{picture}(28,8)
\put(5,3){\circle{10}}
\put(25,3){\circle{10}}
\put(8.5,0){\vector(1,0){12}}
\put(20.5,6){\vector(-1,0){12}}
\put(10,2.5){\vector(-1,0){10}}
\put(2.5,4){\tiny $\epsilon$}
\put(13,8){\tiny $+$}
\put(13,-5){\tiny $+$}
\end{picture}~
}
\newcommand{\fonepPM}{
\,\begin{picture}(28,8)
\put(5,3){\circle{10}}
\put(25,3){\circle{10}}
\put(8.5,0){\vector(1,0){12}}
\put(20.5,6){\vector(-1,0){12}}
\put(10,2.5){\vector(-1,0){10}}
\put(2.5,4){\tiny $\epsilon$}
\put(13,8){\tiny $+$}
\put(13,-5){\tiny $-$}
\end{picture}~
}
\newcommand{\fonepMP}{
\,\begin{picture}(28,8)
\put(5,3){\circle{10}}
\put(25,3){\circle{10}}
\put(8.5,0){\vector(1,0){12}}
\put(20.5,6){\vector(-1,0){12}}
\put(10,2.5){\vector(-1,0){10}}
\put(2.5,4){\tiny $\epsilon$}
\put(13,8){\tiny $-$}
\put(13,-5){\tiny $+$}
\end{picture}~
}
\newcommand{\fonepMM}{
\,\begin{picture}(28,8)
\put(5,3){\circle{10}}
\put(25,3){\circle{10}}
\put(8.5,0){\vector(1,0){12}}
\put(20.5,6){\vector(-1,0){12}}
\put(10,2.5){\vector(-1,0){10}}
\put(2.5,4){\tiny $\epsilon$}
\put(13,8){\tiny $-$}
\put(13,-5){\tiny $-$}
\end{picture}~
}

\newcommand{\ftwop}{
\,\begin{picture}(28,8)
\put(5,3){\circle{10}}
\put(25,3){\circle{10}}
\put(8.5,6){\vector(1,0){12}}
\put(20.5,0){\vector(-1,0){12}}
\put(10,2.5){\vector(-1,0){10}}
\put(2,3.8){\tiny $+$}
\end{picture}~
}

\newcommand{\fonem}{
\,\begin{picture}(28,8)
\put(5,3){\circle{10}}
\put(25,3){\circle{10}}
\put(8.5,0){\vector(1,0){12}}
\put(20.5,6){\vector(-1,0){12}}
\put(10,2.5){\vector(-1,0){10}}
\put(2,3.8){\tiny $-$}
\end{picture}~
}
\newcommand{\ftwom}{
\,\begin{picture}(28,8)
\put(5,3){\circle{10}}
\put(25,3){\circle{10}}
\put(8.5,6){\vector(1,0){12}}
\put(20.5,0){\vector(-1,0){12}}
\put(10,2.5){\vector(-1,0){10}}
\put(2,3.8){\tiny $-$}
\end{picture}~
}

\begin{document}
\title[Vassiliev invariants of 2-bouquet graphs]{Gauss diagram formulas of Vassiliev invariants of spatial 2-bouquet graphs}
\author{Noboru Ito}
\address{
National Institute of Technology, Ibaraki College, 866 Nakane Hitachinaka, Ibaraki 312-8508, Japan
}
\email{nito@ibaraki-ct.ac.jp}
\author{Natsumi Oyamaguchi}
\address{Department of teacher education, 1-1 Daigaku-cho, Yachiyo City, Shumei University, Chiba 276-0003, Japan}
\email{p-oyamaguchi@mailg.shumei-u.ac.jp}
\date{July~9, 2020}
\maketitle
\begin{abstract}
We introduce new formulas that are Vassiliev invariants of flat vertex isotopy classes of spatial 2-bouquet graphs, which are equivalent to 2-string links.  
Although any Gauss diagram formula of Vassiliev invariants of spatial 2-bouquet graphs in a 3-space has been unknown,  
this paper gives the first and simple  example.   
\end{abstract}

\section{Introduction}
A spatial graph is a graph embedded in $\mathbb{R}^3$.  It often becomes a model of molecule as an embedding of a molecular graph, or a coordination polymer (e.g. \cite[Section~1]{BB}).  In general, the interaction between topological graph theory and the investigation of chemical structures is a rich area.   In particular, we would like to emphasize the following two points:
\begin{itemize}
\item Multicyclic polymers having shapes corresponding to rigid $4$-valent graphs (e.g. flat vertex isotopy classes of spatial 2-bouquet graphs) are synthesized \cite{Tezuka2011, Tezuka2015}.   
\item The difference between  spatial graphs affects a condensed matter, e.g. it is shown that the dominance of the trefoil knot in the case of large excluded volumes 
 \cite{UeharaDeguchi2017}.     
\end{itemize}

On the other hand, Deguchi applied the Vassiliev invariant of order two to computational science to study random knotting or linking (1994, e.g. \cite{D, DT}).   One of his motivation, in practical application of the Jones polynomial, is to solve two problems: (1)  Divergence occurs when we evaluate polynomials; (2)  Computational time is growing exponentially with respect to the number of crossings of link diagrams.  

Deguchi \cite{D} gave a solution of two problems showing an algorithm by using the expansion at $q=1$ for the Jones polynomial $V_{K} (q)$: 
\[
V_K (q) = 1 + v_2 (K) \epsilon^2 + v_3 (K) \epsilon^3 +  \dots , \]
where $v_i (K)$ is called Vassiliev invariant of degree $i$.  
We stand on the viewpoint.   
It is also meaningful because 
in general, if all Vassiliev invariants for two knots coincide, then their (Alexander, Conway, Jones, Kauffman, HOMFLY-PT, etc.)  polynomial invariants coincide.  
Nowadays, it is known that every Vassiliev invariant is expressed by a Gauss diagram formula (2000, \cite{gpv}).  It seems likely that this type of formulas is the simplest for computation purposes.   

\begin{figure}[h!]
\includegraphics[width=10cm]{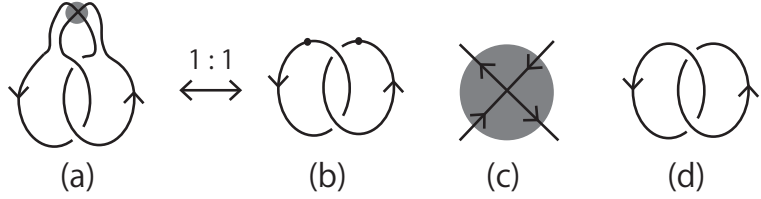}
\caption{(a) : An  oriented bouquet graph
(b) : A $2$-string link with base points ($=$ a $2$-string tangle)
(c) : A neighborhood of the flat vertex
(d) : A link obtained by ignoring the base points of (b).  
}\label{2strings}
\end{figure}

\begin{figure}[h!]
\includegraphics[width=10cm]{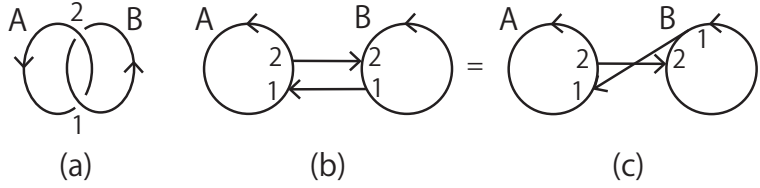}
\caption{(a) : A $2$-component link having components $A$ and $B$, (b) : A Gauss diagram of (a), (c) : Another Gauss diagram of (b).}\label{linkGauss}
\end{figure}

In this paper, we devote ourselves to 2-component link invariants since it is known that there exists one to one correspondence between flat vertex isotopy classes of bouquet graphs to 2-string links (e.g., Figure~\ref{2strings}, for the details of the definitions of bouquet graphs and flat vertex isotopy, see \cite[Section~2]{oyamaguchi}).   
In order to give the statement of main results (Theorems~\ref{thmUS} and \ref{thmOST}), we use definitions of Gauss diagrams and arrow diagrams for links.  For these definitions, please see \cite{ostlund}.  Here, we give an example by Figure~\ref{linkGauss}.    
In the following statements, an estimation ``degree $\le n$"   induced by $n$ arrows is called \emph{order} $n$.    


\begin{theorem}\label{thmUS}
Each of $\langle  \fonep, \cdot \rangle$, $\langle  \ftwop, \cdot \rangle$, $\langle  \fonem, \cdot \rangle$, and $\langle  \ftwom, \cdot \rangle$ is an integer-valued nonzero function that is an invariant of order three of two-component links. 

As a corollary, each of them is also an invariant of order three of spatial graphs in $\mathbb{R}^3$ up to flat vertex isotopy.    
\end{theorem}
\begin{theorem}\label{thmOST}
Each of $\langle  \foneA, \cdot \rangle$, $\langle  \foneB, \cdot \rangle$, and $\langle  \foneC, \cdot \rangle$ $-$ $\frac{1}{3} \langle  \foneD, \cdot \rangle$ is an integer-valued nonzero function that is an invariant of order three of two-component links. 

As a corollary, each of them is also an invariant of order three of spatial graphs in $\mathbb{R}^3$ up to flat vertex isotopy.   
\end{theorem}
\begin{corollary}[\"{O}stlund-Polyak-Viro formula]\label{corOST}
\"{O}stlund-Polyak-Viro formula $\langle  T, \cdot \rangle$, which is
$\langle  \foneA, \cdot \rangle$ $+$ $\langle  \foneB, \cdot \rangle$ $+$ $\langle  \foneC, \cdot \rangle$ $-$ $\frac{1}{3} \langle  \foneD, \cdot \rangle$, becomes a link invariant of order three of two-component links.  
\end{corollary}
Theorem~\ref{USvsOST} implies that invariants in Theorems~\ref{thmUS} and \ref{thmOST} are strictly stronger than \"{O}stlund-Polyak-Viro  formula (Corollary~\ref{corOST}), which is known as a Vassiliev link invariant of degree three \cite{PV, ostlund} (the formula in \cite{PV} is misprinted and \cite{ostlund} gives the correct formula).    
\begin{theorem}\label{USvsOST}
There exists an infinitely many pairs $(i, j)$ of $2$-component links $L_i, L_j$ ($i \neq j$) such that

$\langle \fonep, L_i \rangle \neq \langle \fonep, L_j \rangle,$ $\langle \fonem, L_i \rangle \neq \langle \fonem, L_j \rangle$, 
$\langle \ftwop, L_i \rangle \neq \langle \ftwop, L_j \rangle$, and $\langle \ftwom, L_i \rangle \neq \langle \ftwom, L_j \rangle$ 
for our invariants as in  Theorem~\ref{thmUS} whereas for any pair $i, j$,   
\[
\langle T, L_i \rangle = \langle T, L_j \rangle 
\]
  on \"{O}stlund-Polyak-Viro formula $\langle T, \cdot \rangle$ as in Corollary~\ref{corOST}.   
\end{theorem}

\section{Preliminaries}
If a reader is familiar with the brackets $\langle \cdot, \cdot  \rangle$ and $(\cdot, \cdot)$ introduced in \cite{gpv} or treated in \cite{ostlund},  the reader can skip this section except for Notation~\ref{notation1}.  

\begin{definition}
Let $L$ be a two-component link.  
For $L$, let $D$ be a link diagram and $G$ a (signed oriented) Gauss diagram, where each sign is a local writhe.  In what follows, every Gauss diagram is signed and oriented.  
A sub-Gauss diagram of $G$ is a Gauss diagram obtained by ignoring some arrows.    
Then, let $\sub(G)$ be the set of  sub-Gauss diagrams of $G$. 
For Gauss diagrams $A$ and $z$, $(\cdot, \cdot )$ is an orthonormal scalar product, i.e. $(A, z) = 1$ if $A=z$ and $0$ otherwise.      
Then, 
$\langle \cdot, \cdot \rangle$ is defined by 
\begin{equation}\label{defBracket}
\langle A, G \rangle = \Sigma_{z \in \sub(G)} \sign(z) (A, z),   
\end{equation}
here $\sign(z)$ is the product of the signs in $z$.   
In general, let $\langle S \rangle$ be a $\mathbb{Q}$-vector space generated by the set $S$ of finitely many Gauss diagrams.  
Let $\langle S \rangle$ be a vector space generated by the Gauss diagram having at most $d$ arrows where $d$ is sufficiently large.  
We extend $\langle \cdot, \cdot \rangle$ 
to $\langle S  \rangle \times \langle S  \rangle$ bilinearly.      
\end{definition}
\begin{notation}
In this paper, every circle of Gauss diagrams is oriented counterclockwise.  
When no confusion is likely arise, we omit an orientation on each circle, e.g. $\foneC$.  
\end{notation}
\begin{fact}[The linking number relation {\cite[Page~451, Theorem~5]{PV}}, {\cite[Section~4.1]{ostlund}}]\label{lkfact}
The linking number $lk(L)$ of a two-component link $L$ is given by 
\[
lk(L) = \langle \One, \cdot \rangle = \langle \Oneb, \cdot \rangle.\]
\end{fact}
In particular, for a Gauss diagram $D$ of a two-component link $L$, the above formula is also represented by
\begin{equation}\label{lkeq}
\sum_{z \in \sub(D)} \sum_{\epsilon=+, -} \sign(z) (\Oneep, z) = \sum_{z \in \sub(D)} \sum_{\epsilon=+, -} \sign(z) (\Onebep, z).  
\end{equation}
\begin{notation}\label{notation1}
Let $\epsilon$ be $+$ or $-$ and we fix the sign.    
Let $\fonepEP$ $=$ $\fonepPP + \fonepPM + \fonepMP + \fonepMM$.
\end{notation}
\begin{notation}[Terminological remark for Reidemeister moves]
We use the minimal generating set $\{ \Omega_{1a}, \Omega_{1b},$ $\Omega_{2a}, \Omega_{3a} \}$ of  Reidemeister moves for oriented link diagrams by Polyak \cite[Theorem~1.1]{P}.   
For them, it is convenient to use \"{O}stlund's notations because \cite[Table~1]{ostlund} includes a version involving two component links,  $\Omega_{2a}$ and $\Omega_{3a}$\cite[Theorem~1.1]{P}, represented by Gauss diagrams.  Concretely, $\Omega_{2a}$ corresponds to $\Omega_{II+-}$, and $\Omega_{3a}$ corresponds to the three presentations $\Omega_{III+-+b}$, 
$\Omega_{III+-+m}$, and $\Omega_{III+-+t}$, depending on connectedness of components, for two component links.   

For Reidemeister moves involving one link component, $\Omega_{1a}$ ($\Omega_{1b}$,~resp.) is denoted by $\Omega_{1-+}$ ($\Omega_{1++}$,~resp.)  
\end{notation}
\begin{figure}[h!]
\includegraphics[width=13cm]{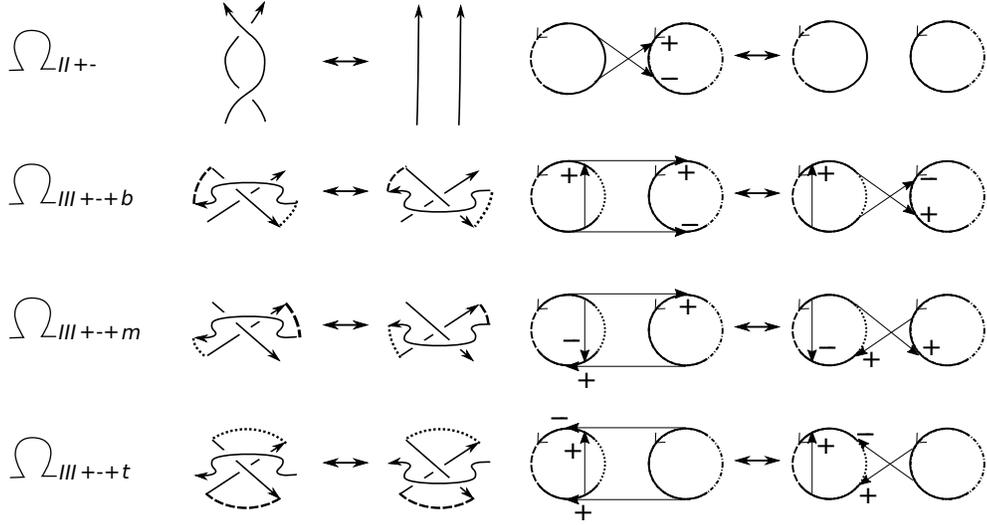}
\caption{Reidemeister moves for 2-components and their Gauss diagrams.  $\Omega_{II +-}$ (upper line), $\Omega_{III+-+b}$, $\Omega_{III+-+m}$, $\Omega_{III+-+t}$ (lower line), which are included in \cite[Table~1]{ostlund}.}\label{reide&relator}
\end{figure}

\section{Proof of main results}
We will show Theorem~\ref{thmUS} after proving Theorem~\ref{thmOST}.   In this section, we freely use \"{O}stlund's notation \cite{ostlund} of \emph{arrow diagrams} and their moves (in particular, cf.~\cite[Section~1.6 and Table~1]{ostlund}) since \cite{ostlund} is the paper giving the proof of all statements of \cite{PV}.   If a  reader is familiar with word-theoretic approach to Gauss diagrams, an advantage for computation is given as in \cite{Ito}.     
Denoted by $D^{\Omega_*}_l$ ($D^{\Omega_*}_r$,~resp.) the left (right,~resp.) Gauss diagram of in each Reidemeister move $\Omega_*$ as in \cite[Table~1]{ostlund}.     
\subsection{Proof of Theorem~\ref{thmOST}}
Since the invariance under each of $\Omega_{1+ \pm}$, $\Omega_{2+ -}$, and $\Omega_{3 + - + \pm}$ is obvious, we show the invariance under $\Omega_{II +-}$, $\Omega_{III +-+ b}$, $\Omega_{III+-+ m}$, and $\Omega_{III+-+ t}$.  

\noindent$\bullet$ $\Omega_{II +-}$.  
First we consider $\foneC$.  
\begin{align*}
\langle &\foneC, D_l^{\Omega_{II+-}}  \rangle - \langle \foneC,   D_r^{\Omega_{II+-}} \rangle\\ 
&= \!\!\!\! \sum_{z^{(l)} \in \sub \left(D_l^{\Omega_{II+-}} \right)} \sign(z^{(l)}) (\foneC, z^{(l)})  - \!\!\!\! \sum_{z^{(r)} \in \sub \left(D_r^{\Omega_{II+-}} \right)} \sign(z^{(r)}) (\foneC, z^{(r)}).  
\end{align*}
Since $D_l^{\Omega_{II+-}}$ has two more arrows than  $D_r^{\Omega_{II+-}}$, $z^{(l)}$ is denoted by $z^{(l)}_2$ if $z^{(l)}$ contains these two arrows.  Similarly, if $z^{(l)}$ contains $i$ ($i=0, 1$) arrow of these two arrows, $z^{(l)}$ is denoted by $z^{(l)}_i$.  Then, 
\begin{align*}
&\sum_{z^{(l)} \in \sub \left(D_l^{\Omega_{II+-}} \right)} \sign(z^{(l)}) (\foneC, z^{(l)}) - \sum_{z^{(r)} \in \sub \left(D_r^{\Omega_{II+-}} \right)} \sign(z^{(r)}) (\foneC, z^{(r)})\\
&= \sum_{i=0}^2 \sum_{z^{(l)}_i \in \sub \left(D_l^{\Omega_{II+-}} \right)} \sign(z^{(l)}_i) (\foneC, z^{(l)}_i) - \!\!\!\!\!\! \sum_{z^{(r)} \in \sub \left(D_r^{\Omega_{II+-}} \right)} \!\!\!\!\!\! \sign(z^{(r)}) (\foneC, z^{(r)}).  
\end{align*}
Here, by definition, the set of elements labeled by $z^{(l)}_0$ corresponds bijectively to that of $z^{(r)}$.   Then, 
\begin{align*}
\sum_{z^{(l)}_0 \in \sub \left(D_l^{\Omega_{II+-}} \right)} \sign(z^{(l)}_0) (\foneC, z^{(l)}_0) = \sum_{z^{(r)} \in \sub \left(D_r^{\Omega_{II+-}} \right)} \sign(z^{(r)}) (\foneC, z^{(r)}).
\end{align*}
In general, for any linear sum $A$ of Gauss diagrams, 
\begin{align}\label{z0Eq}
\sum_{z^{(l)}_0 \in \sub \left(D_l^{\Omega_{II+-}} \right)} \sign(z^{(l)}_0) (A, z^{(l)}_0) = \sum_{z^{(r)} \in \sub \left(D_r^{\Omega_{II+-}} \right)} \sign(z^{(r)}) (A, z^{(r)}).
\end{align}
Note also that two arrows in the difference between $D_l^{\Omega_{II+-}}$ and $D_r^{\Omega_{II+-}}$ have $+$ and $-$ signs, respectively.  Then if $z^{(l)}_1$ includes $+$ sign in the difference between $D_l^{\Omega_{II+-}}$ and $D_r^{\Omega_{II+-}}$, we denote it by $z^+_1$.  
Then, there exists another $z^{(l)}_1$ with $-$ sign in the difference, and we denote it by $z^-_1$.  Then, 
\begin{equation*}
\begin{split}
&\sum_{z^{(l)}_1 \in \sub \left(D_l^{\Omega_{II+-}} \right)} \sign(z^{(l)}_1) (\foneC, z^{(l)}_1) \\
&= \sum_{z^+_1 \in \sub \left(D_l^{\Omega_{II+-}} \right)} \sign(z^+_1) (\foneC, z^+_1) + \sum_{z^-_1 \in \sub \left(D_l^{\Omega_{II+-}} \right)} \sign(z^-_1) (\foneC, z^-_1)\\
&= \sum_{z^+_1 \in \sub \left(D_l^{\Omega_{II+-}} \right)} \sign(z^+_1) (\foneC, z^+_1) -  \sum_{z^+_1 \in \sub \left(D_l^{\Omega_{II+-}} \right)} \sign(z^+_1) (\foneC, z^+_1)\\
&= 0.  
\end{split}
\end{equation*}
In general, for any linear sum $A$ of Gauss diagrams, 
\begin{equation}\label{z1Eq}
\begin{split}
&\sum_{z^{(l)}_1 \in \sub \left(D_l^{\Omega_{II+-}} \right)} \sign(z^{(l)}_1) (A, z^{(l)}_1) 
= 0.  
\end{split}
\end{equation}

Next, by (\ref{z0Eq}), 
\begin{align*}
\langle &\foneC-\frac{1}{3} \foneD, D_l^{\Omega_{II+-}}  \rangle - \langle \foneC-\frac{1}{3} \foneD,   D_r^{\Omega_{II+-}} \rangle\\ 
&= \langle \foneC, D_l^{\Omega_{II+-}} \rangle - \frac{1}{3} \langle  \foneD, D_l^{\Omega_{II+-}} \rangle. 
\end{align*}
Further, by (\ref{z1Eq}), 
\begin{align*}
&\langle \foneC, D_l^{\Omega_{II+-}} \rangle - \frac{1}{3} \langle  \foneD, D_l^{\Omega_{II+-}} \rangle\\
&=  \sum_{z^{(l)}_2 \in \sub \left(D_l^{\Omega_{II+-}} \right)} \sum_{\epsilon=+, -}  \sign (z^{(l)}_2) (\foneCII, z^{(l)}_2)  
\\ & \qquad \qquad \qquad \qquad \qquad \qquad 
- 
\sum_{\tilde{z}^{(l)}_2 \in \sub \left(D_l^{\Omega_{II+-}} \right)} \sum_{\epsilon=+, -} 
 \sign (\tilde{z}^{(l)}_2) (\foneDII, \tilde{z}^{(l)}_2)  \\
& ({\textrm{for the second term, $\frac{1}{3} \cdot 3$ appears as in \cite[Sec.~4.8.3]{ostlund}
by symmetry of }} \foneD)  \\
&= 0 \qquad (\because  {\textrm{The linking number relation~}}(\ref{lkeq}) ).
\end{align*}
We note also that by (\ref{z0Eq}) and (\ref{z1Eq}), it is easy to see that the following two formulas  hold
\[\langle \foneA, D_l^{\Omega_{II+-}} -  D_r^{\Omega_{II+-}} \rangle=0,
\]
and
\[
\langle \foneB, D_l^{\Omega_{II+-}} -  D_r^{\Omega_{II+-}} \rangle =0.
\]

Below, we discuss $\Omega_{III+-+ *}$ ($*=b, m, t$).   The difference between $D_l^{\Omega_{III+-+*}}$ and $D_r^{\Omega_{III+-+*}}$ is the addition of three  arrows.  Let $z^{(l)}$ ($z^{(r)}$,~resp.) be an element of $ \sub(D_l^{\Omega_{III+-+*}})$ ($\sub(D_r^{\Omega_{III+-+*}})$,~resp.).  
Then, if $z^{(l)}$ ($z^{(r)}$,~resp.)  contains $i$ ($i=0, 1, 2, 3$) arrow(s) of these three arrows, $z^{(l)}$ ($z^{(r)}$,~resp.) is denoted by $z^{(l)}_i$ ($z^{(r)}_i$,~resp.).

By the same argument as the case $\Omega_{II+-}$, for $i=0, 1$, we have
\begin{align}\label{z0EqIII}
\sum_{z^{(l)}_0 \in \sub \left(D_l^{\Omega_{III+-+*}} \right)} \sign(z^{(l)}_0) (A, z^{(l)}_0) = \sum_{z^{(r)}_0 \in \sub \left(D_r^{\Omega_{III+-+*}} \right)} \sign(z^{(r)}_0) (A, z^{(r)}_0).
\end{align}
and 
\begin{align}\label{z1EqIII}
\sum_{z^{(l)}_1 \in \sub \left(D_l^{\Omega_{III+-+*}} \right)} \sign(z^{(l)}_1) (A, z^{(l)}_1) = \sum_{z^{(r)}_1 \in \sub \left(D_r^{\Omega_{III+-+*}} \right)} \sign(z^{(r)}_1) (A, z^{(r)}_1).
\end{align}
Then, we note that by (\ref{z0EqIII}) and (\ref{z1EqIII}), the following two formulas hold.  For any $*$ ($=b, m, t$), 
\[\langle \foneA, D_l^{\Omega_{III+-+*}} -  D_r^{\Omega_{II+-+*}} \rangle=0,
\]
and
\[
\langle \foneB, D_l^{\Omega_{III+-+*}} -  D_r^{\Omega_{III+-+*}} \rangle =0.
\]
Hence, we discuss $\foneC$ $-\frac{1}{3} \foneD$ in the following.  

\noindent$\bullet$ $\Omega_{III+-+ b}$.

By (\ref{z1Eq}), 
\begin{align*}
&\langle \foneC, D_l^{\Omega_{III+-+b}} \rangle - \frac{1}{3} \langle  \foneD, D_l^{\Omega_{III+-+b}} \rangle\\
&=  \sum_{z^{(l)}_2 \in \sub \left(D_l^{\Omega_{III+-+b}} \right)} \sum_{\epsilon=+, -}  \sign (z^{(l)}_2) (\foneCII, z^{(l)}_2)  
\\ & \qquad \qquad \qquad \qquad \qquad \qquad 
- 
\sum_{\tilde{z}^{(l)}_2 \in \sub \left(D_l^{\Omega_{III+-+b}} \right)} \sum_{\epsilon=+, -} 
 \sign (\tilde{z}^{(l)}_2) (\foneDII, \tilde{z}^{(l)}_2)  \\
& ({\textrm{for the second term, $\frac{1}{3} \cdot 3$ appears as in \cite[Sec.~4.8.3]{ostlund}
by symmetry of }} \foneD)  \\
&= 0 \qquad (\because  {\textrm{The linking number relation~}}(\ref{lkeq}) ).
\end{align*}

\noindent$\bullet$ $\Omega_{III+-+ m}$.  Since there is no element corresponding to $z^{(l)}_2$ or $z^{(r)}_2$, 
\[
\langle \foneC, D_l^{\Omega_{III+-+m}} -  D_r^{\Omega_{III+-+m}} \rangle =0, 
\]
and 
\[
\langle \foneD, D_l^{\Omega_{III+-+m}} -  D_r^{\Omega_{III+-+m}} \rangle =0.
\]

\noindent$\bullet$ $\Omega_{III+-+ t}$.  

Since we have the same formulas  as in the above $\Omega_{III+-+ b}$ case  except for replacing $b$ with $t$, we omit its proof.  
$\hfill\Box$
\subsection{Proof of Theorem~\ref{thmUS}}
Since the invariance under each of $\Omega_{1+ \pm}$, $\Omega_{2+ -}$, and $\Omega_{3 + - + \pm}$ is obvious, we show the invariance under $\Omega_{II +-}$, $\Omega_{III +-+ b}$, $\Omega_{III+-+ m}$, and $\Omega_{III+-+ t}$.  Below, the same notations $z^{(l)}_i$, $z^{(r)}_i$, and $z^{\pm}_1$ as those of the proof of Theorem~\ref{thmOST} apply.  

\noindent$\bullet$ $\Omega_{II +-}$.  

\begin{align*}
\langle &\fonep, D_l^{\Omega_{II+-}}  \rangle - \langle \fonep,   D_r^{\Omega_{II+-}} \rangle\\ 
&= \!\!\!\! \sum_{z^{(l)} \in \sub \left(D_l^{\Omega_{II+-}} \right)} \sign(z^{(l)}) (\fonep, z^{(l)})  - \!\!\!\! \sum_{z^{(r)} \in \sub \left(D_r^{\Omega_{II+-}} \right)} \sign(z^{(r)}) (\fonep, z^{(r)})\\ 
&= \sum_{i=0}^2 \sum_{z^{(l)}_i \in \sub \left(D_l^{\Omega_{II+-}} \right)} \!\!\! \sign(z^{(l)}_i) (\fonep, z^{(l)}_i) - \!\!\!\!\!\!\!\!\!\!  \sum_{z^{(r)} \in \sub \left(D_r^{\Omega_{II+-}} \right)} \sign(z^{(r)}) (\fonep, z^{(r)}) \\
&= \sum_{z^{(l)}_2 \in \sub \left(D_l^{\Omega_{II+-}} \right)} \!\!\! \sign(z^{(l)}_2) (\fonep, z^{(l)}_2)   
\qquad (\because (\ref{z0Eq}), (\ref{z1Eq})) \\
&= 0 \quad (\because~{\text{there is no $z_2$ that takes non-zero value}}).  
\end{align*}

\begin{figure}[h!]
\includegraphics[width=12.5cm]{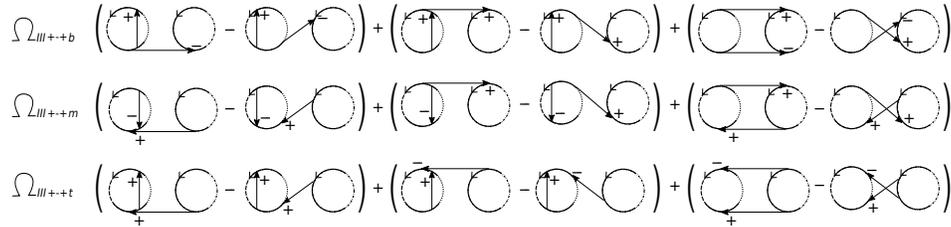}
\caption{Three pairings corresponding to three subtractions for each Reidemeister move $\Omega_{III+-+ *}$.}\label{6term}
\end{figure}

\begin{figure}[h!]
\includegraphics[width=12cm]{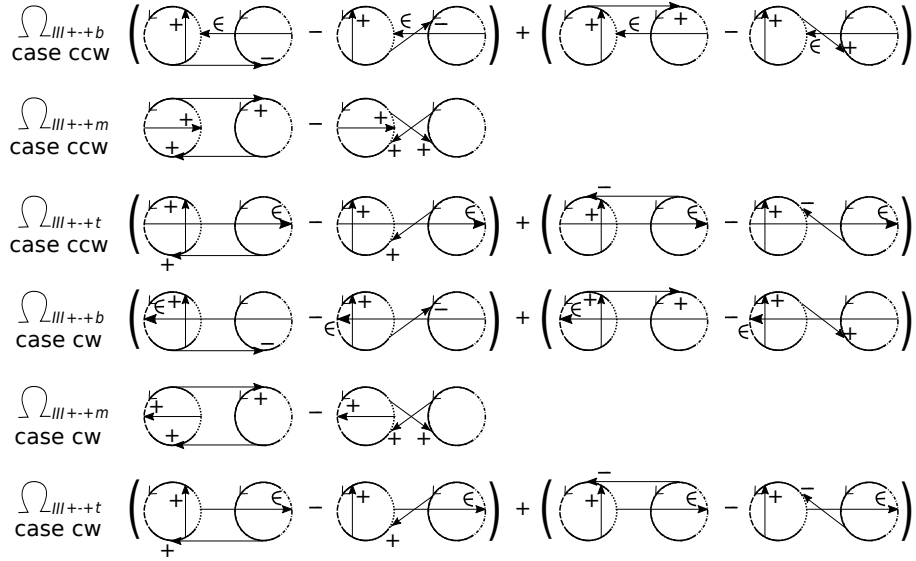}
\caption{This is a list of the possible terms that take non-trivial values in the form as in Figure~\ref{6term}, where $\epsilon$ $=$ $+$ or $-$.}\label{6term2}
\end{figure}

Next, we consider $\Omega_{III+-+ *}$ ($*=b, m, t$). 
We note that $z^{(l)}_2$ corresponds to $z^{(r)}_2$ under the Reidemeister move $\Omega_{III+-+*}$ ($*=b, m, t$).  Then, each $\Omega_{III+-+*}$ corresponds to a sum of the three canonical subtractions ``$z^{(l)}_2 - z^{(r)}_2$'' of pairings $(z^{(l)}_2, z^{(r)}_2)$ as in Figure~\ref{6term}.
Further, by the definition of $(\fonep, \cdot)$, each term, which survives in subtractions ``$z^{(l)}_2$ $-$ $z^{(r)}_2$",  consists of two arrows relevant to $\Omega_{III+-+*}$ and the other one with a sign $\epsilon$ (for $*=b, t$) or $+$ (for $*=m$) as in the case~``ccw" (i.e., ``{\bf c}ounter{\bf c}lock{\bf w}ise as in $\Two$" embedded in $\fonep$) of Figure~\ref{6term2}. 


\noindent $\bullet$ $\Omega_{III+-+ b}$.  
\begin{align*}
\langle &\fonep, D_l^{\Omega_{III+-+b}}  \rangle - \langle \fonep,   D_r^{\Omega_{III+-+b}} \rangle\\ 
&= \!\!\!\! \sum_{z^{(l)} \in \sub \left(D_l^{\Omega_{III+-+b}} \right)} \!\!\!\!\!\! \sign(z^{(l)}) (\fonep, z^{(l)})  - \!\!\!\! \sum_{z^{(r)} \in \sub \left(D_r^{\Omega_{III+-+b}} \right)}  \!\!\!\!\!\! \sign(z^{(r)}) (\fonep, z^{(r)})\\ 
 \end{align*}
\begin{align*}
&= \sum_{i=0}^3 \left(\sum_{z^{(l)}_i \in \sub \left(D_l^{\Omega_{III+-+b}} \right)} \!\!\!\!\!\!\!\!\!\!\!\!\! \sign(z^{(l)}_i) (\fonep, z^{(l)}_i) - \!\!\!\!\!\!\!\!\!\!  \sum_{z_i^{(r)} \in \sub \left(D_r^{\Omega_{III+-+b}} \right)} \!\!\!\!\!\!\!\!\!\!\!\!\! \sign(z_i^{(r)}) (\fonep, z_i^{(r)}) \right) \\
&= \sum_{z^{(l)}_2 \in \sub \left(D_l^{\Omega_{III+-+b}} \right)} \!\!\!\!\!\! \sign(z^{(l)}_2) (\fonep, z^{(l)}_2) - \!\!\!\!\!\!\!\!\!\!  \sum_{z_2^{(r)} \in \sub \left(D_r^{\Omega_{III+-+b}} \right)} \!\!\!\!\!\! \sign(z_2^{(r)}) (\fonep, z_2^{(r)})  \\
&\qquad (\because (\ref{z0Eq}), (\ref{z1Eq}),~{\text{there is no $z_3$ that takes non-zero value}})\\
&=\left( \fonep, \sum_{z^{(l)}_2 \in \sub \left(D_l^{\Omega_{III+-+b}} \right)}
 \sign(z^{(l)}_2) z^{(l)}_2  - \sum_{z^{(r)}_2 \in \sub \left(D_r^{\Omega_{III+-+b}} \right)}
 \sign(z^{(r)}_2) z^{(r)}_2 \right)\\
 &= \left( \fonep, \sum_{{\text{pairings}}~(z^{(l)}_2,~ z^{(r)}_2)~{\text{as in Figure~\ref{6term} }}} \left(
 \sign(z^{(l)}_2) z^{(l)}_2  - 
 \sign(z^{(r)}_2) z^{(r)}_2 \right) \right)\\ 
&= \left( \fonep, \sum_{z^{(l)}_2,~ \tilde{z}^{(l)}_2 {\text{surviving on $\Omega_{+-+b}$ as in Figure~\ref{6term2} }}} \left(
 \sign(z^{(l)}_2) z^{(l)}_2 + \sign(\tilde{z}^{(l)}_2) \tilde{z}^{(l)}_2 \right) \right) \\
&= \sum_{z^{(l)}_2,~ \tilde{z}^{(l)}_2 {\text{ surviving on $\Omega_{+-+b}$ as in Figure~\ref{6term2} }}} \left(
 \epsilon ( \fonep, z^{(l)}_2)  - \epsilon ( \fonep, \tilde{z}^{(l)}_2 ) \right)\\
 &= 0.      
\end{align*}

In the same way as the above, using Figure~\ref{6term} and the ``case~ccw" as in Figure~\ref{6term2}, 
$\langle \fonep, D_l^{\Omega_{III+-+m}}  \rangle$ $-$ $\langle \fonep,   D_r^{\Omega_{III+-+m}} \rangle$ $=$ $0$ and  $\langle \fonep, D_l^{\Omega_{III+-+t}}  \rangle$ $-$ $\langle \fonep,   D_r^{\Omega_{III+-+t}} \rangle$ $=$ $0$.   
Hence, $\fonep$ is a link invariant.  
This fact together with  Theorem~\ref{thmOST} implies that $\fonem$ ($=$ $\fone-\fonep$) is also a link invariant.  

By Figure~\ref{6term} and ``case~cw" (i.e. clockwise) as in Figure~\ref{6term2}, the same argument of the proof of the invariance of $\fonep$ applies, we have a proof of case $\ftwop$. 
The invariance of $\ftwop$ implies that $\ftwom$ ($=$ $\foneB-\ftwop$) is also a link invariant.
$\hfill\Box$ 
\section{Proof of Theorem~\ref{USvsOST}}
\begin{figure}[h!]
\includegraphics[width=12.5cm]{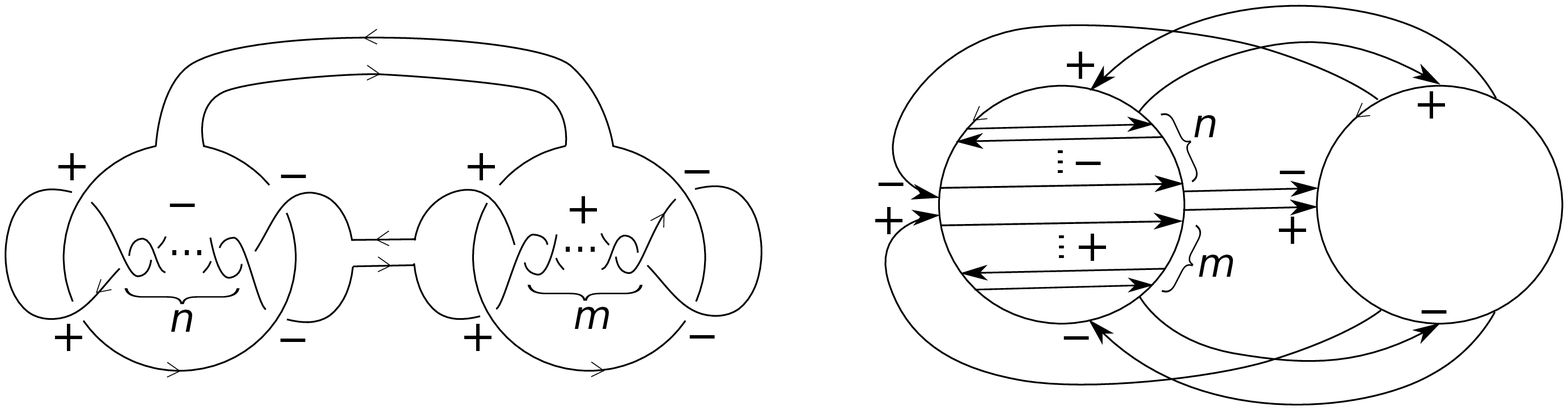}
\caption{For $L(m, n)$ with a fixed pair $(m, n)$, the left figure is a link diagram and the right figure is its Gauss diagram; $m$ ($n$,~resp.) denotes $m$ ($n$,~resp.) crossings.}\label{ex2}
\end{figure}

\begin{figure}[h!]
\includegraphics[width=12.5cm]{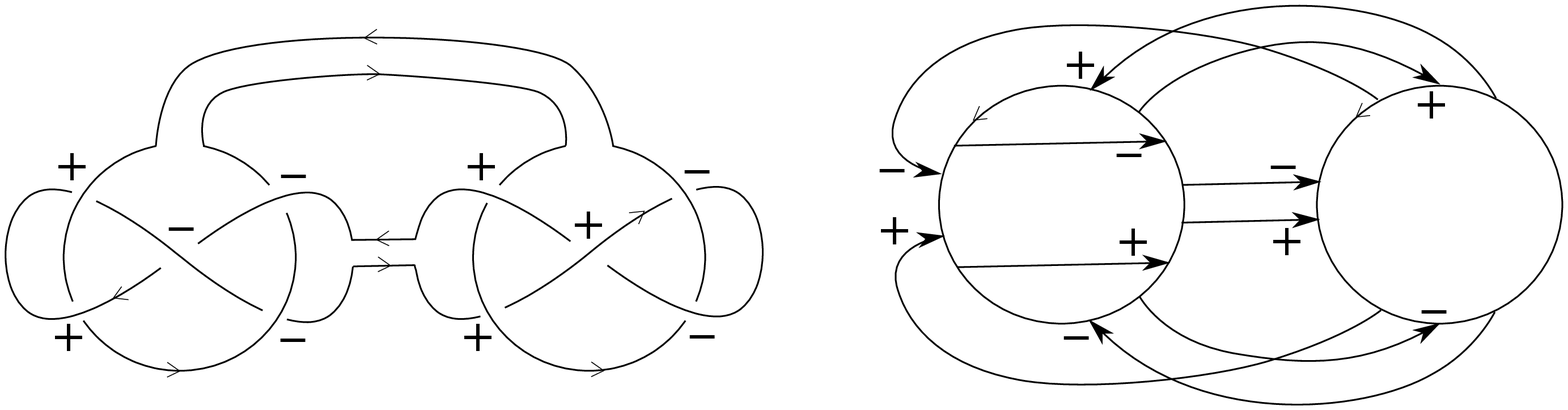}
\caption{$L(1, 1)$.}\label{ex1}
\end{figure}

Let $m$ and $n$ be odd positive integers.  Let $L(m, n)$ be a $2$-component link with $m+n+8$ as in Figure~\ref{ex2} (e.g. for $L(1, 1)$, see Figure~\ref{ex1}).  

Let $L_n$ $=$ $L(n, n)$.  By Table~\ref{tableCompared}, for $L_n$, the \"{O}stlund-Polyak-Viro  formula (Corollary~\ref{corOST}) takes the same value ($=0$) as that of $L_m$ even if $m \neq n$  whereas our invariants take values $-n$ or $n$, which implies the statement of Theorem~\ref{USvsOST}.   
$\hfill\Box$
 
\begin{table}[h!]
\caption{}\label{tableCompared}
\[
\begin{tabular}{|c|c|c|c|c|}\hline
Invariants of Theorem~\ref{thmUS} &  \fonep & \ftwop & \fonem & \ftwom   \\ \hline
Values of $L(m, n)$&$-n$&$-n$&$m$&$m$ \\ \hline
\end{tabular}
\]
\[
\begin{tabular}{|c|c|c|c|}\hline
Invariants of Theorem~\ref{thmOST}    &  \foneA & \foneB & $\foneC - \frac{1}{3} \foneD$   \\ \hline
Values of $L(m, n)$&$m-n$&$m-n$&$0$ \\ \hline
\end{tabular}
\]
\[
\begin{tabular}{|c|c|}\hline
\"{O}stlund-Polyak-Viro formula & \foneA + \foneB + $\foneC - \frac{1}{3} \foneD$  \\  
of Corollary~\ref{corOST} & \\ \hline
Values of $L(m, n)$& $2(m-n)$ \\ \hline
\end{tabular}
\]
\end{table}  

\section*{Acknowledgements}
We would like to thank Professor  Tetsuo Deguchi and Professor  Erica  Uehara for their comments.  We also thank Professor Koya Shimokawa for his guidance and comments.

\end{document}